\documentclass[11pt]{article}
\usepackage{graphicx, colordvi}
\usepackage{amsfonts}
\usepackage{amssymb}
\usepackage{amsthm,cite,color}
\usepackage{dsfont}
\usepackage{epsfig}
\usepackage{mathrsfs}
\usepackage{amsfonts}
\usepackage{amssymb}
\usepackage{amsmath}
\usepackage{amssymb,amsfonts,amsmath,amsthm,cite,color}
\usepackage{dsfont}
\usepackage{epsfig}
\usepackage{mathrsfs}
\usepackage{longtable}

\newtheorem{theo}{Theorem}
\newtheorem{rem}{Remark}

\newtheorem{coro}[theo]{Corollary}

\newtheorem{lem}[theo]{Lemma}

\newtheorem{coj}[theo]{Conjecture}
\makeatletter \@addtoreset{equation}{section}
\@addtoreset{theo}{section} \makeatother

\newcommand{\bN} { {\mathbb{N}}}

\newcommand{\bZ} { {\mathbb{Z}}}

\newcommand{\la} { {\langle}}
\newcommand{\ra} { {\rangle}}

\def\qed{\hfill \rule{4pt}{7pt}}
\def\pf{\noindent {\it Proof.} }

\begin{document}
\begin{center}

 {\Large \bf Congruences Related to Dual Sequences\\ [10pt] and Catalan Numbers
}
\end{center}
\begin{center}
{  Rong-Hua Wang}$^{1}$ and {Michael X.X. Zhong}$^{2}$

   $^1$School of Mathematical Sciences\\
   Tiangong University \\
   Tianjin 300387, P.R. China \\
   wangronghua@tiangong.edu.cn \\[10pt]

   $^2$School of Science\\
   Tianjin University of Technology \\
   Tianjin 300384, P.R. China\\
   zhong.m@tjut.edu.cn

\end{center}

\vskip 6mm \noindent {\bf Abstract.}
During the study of dual sequences, Sun introduced the polynomials
\[
D_n(x,y)=\sum_{k=0}^{n}{n\choose k}{x\choose k}y^k\text{ and }
S_n(x,y)=\sum_{k=0}^{n}\binom{n}{k}\binom{x}{k}\binom{-1-x}{k} y^k.
\]
Many related congruences have been established and conjectured by Sun.
Here we generalize some of them by determining
\[
\sum_{k=0}^{p-1}D_k(x_1,y_1)D_k(x_2,y_2)\pmod p
\text{ and }
\sum_{k=0}^{p-1}S_k(x_1,y_1)S_k(x_2,y_2)\pmod p
\]
for any odd prime $p$ and $p$-adic integers $x_i,\ y_i$ with $i\in\{1,2\}$.
Considering the immediate connection between binomial coefficients and Catalan numbers, we also characterize
\[
\sum_{n=0}^{p-1}\left(\sum_{k=0}^n {n \choose k} \frac{C_k}{a^k}\right)^2
\pmod {p},
\]
where $C_k$ denotes the $k$th Catalan number, $a\in\mathbb{Z}\setminus \{0\}$ with $\gcd(a,p)=1$.
These confirm and generalise some of Sun's conjectures.

\noindent {\bf Keywords}: Congruence; Dual sequence; Catalan number.

\section{Introduction}
Let $\{a_n\}$ be a sequence of numbers,
its \emph{dual sequence} $\{a_n^*\}$ is given as
\[
a_n^*=\sum_{k=0}^n{n \choose k}(-1)^k a_k
\]
for $n\in\mathbb{N}$.
One can consult \cite{GKP1989, Sun2003,KW2006,Sun201512} for properties of dual sequences and related combinatorial identities.
Let $a_n=(-1)^n\binom{x}{n}$, it is easy to see that
\[
a_n^*=\sum_{k=0}^n{n \choose k}\binom{x}{k}=\binom{n+x}{n}=(-1)^n\binom{-1-x}{n}.
\]

In \cite{RV2003}, Rodriguez-Villegas presented four conjectures on
\[
\sum_{k=0}^{p-1}\binom{x}{k}\binom{-1-x}{k} \pmod {p^2}
\]
for $x\in\{-\frac{1}{2},-\frac{1}{3},-\frac{1}{4},-\frac{1}{6}\}$ and any prime $p>3$.
These conjectures were proved by Mortenson \cite{Mortenson2003a,Mortenson2003b}.

In 2006, during the study of special values of spectral zeta functions, Kimoto and Wakayama \cite{KW2006} introduced $\tilde{J}_2(n)$ as an analogue of the Ap\'ery numbers, which can be expressed as
\[
\tilde{J}_2(n)=\sum_{k=0}^n{n\choose k}(-1)^k{-1/2\choose k}^2.
\]
They studied the congruence properties of $\tilde{J}_2(n)$ similar to those of the Ap\'ery numbers
and conjectured that
\begin{equation*}
\sum_{n=0}^{p-1}\tilde{J}_2(n)^2\equiv \left(\frac{-1}{p}\right) \pmod {p^3},
\end{equation*}
where $\left(\frac{\bullet}{p}\right)$ is the Legendre symbol.
This conjecture was confirmed by Long, Osburn and Swisher \cite{LOS2014} in 2014.

Inspired by these works, Sun~\cite{Sun201512} introduced the polynomials
\[
D_n(x,y)=\sum_{k=0}^{n}{n\choose k}{x\choose k}y^k \text{ and }
S_n(x,y)=\sum_{k=0}^{n}\binom{n}{k}\binom{x}{k}\binom{-1-x}{k} y^k.
\]
%
Many congruences involving $D_n(x,y)$ and $S_n(x,y)$ were established in \cite{Sun201512}.
Specifically, Sun showed that
\begin{equation}\label{EQ:D_n(x,y)}
\sum_{k=0}^{p-1}D_k(x,y)D_k(-1-x,y)\equiv (-1)^{\la x \ra_p} \pmod{p}
\end{equation}
and
\begin{align}\label{EQ:S_n(x,y)}
&\sum_{k=0}^{p-1}S_k(x,y)^2
\equiv
\left\{
  \begin{array}{ll}
   \left(\frac{-1}{p}\right) \pmod p & \hbox{ if } x \equiv -1/2 \pmod p,
   \\[10pt]
   0 \pmod p  & \hbox{ otherwise, }
  \end{array}
\right.
\end{align}
where $p$ is an odd prime, $x,y$ are $p$-adic integers with $y\not\equiv 0\pmod {p}$ and $\la x\ra_p$ denotes the least $r\in\bN$ such that $x\equiv r\pmod{p}$.
Besides, quite a few conjectures on congruence relations and integer-valuedness involving $D_n(x,y)$ or $S_n(x,y)$ were posed.
During the past few years, many of them have been proved,
see \cite{Guo2016, HouWang2018, Liu2016, Liu2017 }, while some of them are still remaining open.

In this paper, we first consider the congruences on $D_n(x,y)$ and $S_n(x,y)$ in a more general situation compared with identities  \eqref{EQ:D_n(x,y)} and \eqref{EQ:S_n(x,y)}.
%
\begin{theo}\label{TH:Main1}
Let $p$ be an odd prime, $x_i,y_i$ be any $p$-adic integers and $r_i=\langle x_i \rangle_p$, where $i\in\{1,2\}$.
Then we have
\begin{equation}\label{EQ:main1}
\sum_{k=0}^{p-1}D_k(x_1,y_1)D_k(x_2,y_2)
\equiv\left\{
  \begin{array}{ll}
   0 \pmod p & \hbox{ if } r<p-1,
   \\[10pt]
   y_1^{r_1}(-y_2)^{r_2} \pmod p  & \hbox{ if } r=p-1,
  \end{array}
\right.
\end{equation}
where $r=r_1+r_2$, and
\begin{align}
 &\sum_{k=0}^{p-1}S_k(x_1,y_1)S_k(x_2,y_2)
\equiv\left\{
  \begin{array}{ll}
   \left(\frac{-y_1y_2}{p}\right) \pmod p & \hbox{ if } r_1\!=\!r_2\!=\!\frac{p-1}{2},\!
   \\[10pt]
   0 \pmod p  & \hbox{ otherwise }.
  \end{array}\label{EQ:main2}
\right.
\end{align}

\end{theo}

It is easy to see that equality~\eqref{EQ:D_n(x,y)} is the special case of identity~\eqref{EQ:main1} when $x_1+x_2=-1$ and $y_1=y_2$,
and equality~\eqref{EQ:S_n(x,y)} follows from identity~\eqref{EQ:main2} by taking $x_1=x_2=x$ and $y_1=y_2=y$.

By taking appropriate special values of $x_i$ and $y_i$, one can get other interesting congruences.
In section \ref{sec:proof1}, we will use Theorem \ref{TH:Main1} in this way to confirm and generalise the following conjecture posed by Sun.
\begin{coj}\cite[Conjecture 6.7]{Sun201512}
Let $p>3$ be a prime. We have
\begin{equation}\label{CONJ:sun6.7}
\sum_{n=0}^{p-1}\left(\sum_{k=0}^n {n \choose k}\frac{{2k \choose k}}{2^k}\right)
                \left(\sum_{k=0}^n {n \choose k}\frac{{2k \choose k}}{(-6)^k}\right)
\equiv \left(\frac{3}{p}\right)
\pmod {p}.
\end{equation}
\end{coj}

By the close relation between Catalan numbers and binomial coefficients,
it is interesting to explore whether there are similar congruence properties of Catalan numbers.
In fact, we derive the following congruence equality.

\begin{theo}\label{TH:Main3}
Let $p$ be a prime and $a$ an integer such that $\gcd(a,p)=1$. Then \begin{equation}\label{EQ:a}
\sum_{n=0}^{p-1}\left(\sum_{k=0}^n {n \choose k} \frac{C_k}{a^k}\right)^2
\equiv 4\left(\frac{-1}{p}\right)-\frac{\delta}{a}\left(\frac{\delta}{p}\right)+a+1 \pmod {p},
\end{equation}
where $C_k$ denotes the $k$th Catalan number, $a\in\mathbb{Z}\setminus \{0\}$, $\delta=a(a+4)$ and $\left(\frac{\bullet}{p}\right)$ is the Legendre symbol.
\end{theo}
In particular, by setting $a=2$ and $a=-6$ respectively in identity~\eqref{EQ:a}, we confirm another conjecture by Sun.
\begin{coj}\cite[Conjecture 6.7]{Sun201512}
Let $p>3$ be a prime. We have
\begin{equation}\label{EQ:2}
\sum_{n=0}^{p-1}\left(\sum_{k=0}^n {n \choose k} \frac{C_k}{2^k}\right)^2
\equiv 4\left(\frac{-1}{p}\right)-6\left(\frac{3}{p}\right)+3 \pmod {p}
\end{equation}
and
\begin{equation}\label{EQ:-6}
\sum_{n=0}^{p-1}\left(\sum_{k=0}^n {n \choose k} \frac{C_k}{(-6)^k}\right)^2
\equiv 4\left(\frac{-1}{p}\right)+2\left(\frac{3}{p}\right)-5 \pmod {p}.
\end{equation}
\end{coj}

\section{Proof of Theorem~\ref{TH:Main1}}\label{sec:proof1}
To prove Theorem~\ref{TH:Main1}, we first need two Lemmas.
\begin{lem}\label{LM:NewFormD}
For any $n\in\mathbb{N}$ we have
\begin{equation*}\label{EQ:D_n}
D_n(x,y)=\sum_{k=0}^{n}{n\choose k}{x+k\choose k}y^k(1-y)^{n-k}.
\end{equation*}
\end{lem}
Lemma~\ref{LM:NewFormD} is a special case of the following identity
\[
\sum_{k=0}^{n}\binom{n}{k}\binom{z+k}{k}(x-y)^{n-k}y^k=
\sum_{k=0}^{n}\binom{n}{k}\binom{z}{k}x^{n-k}y^k,
\]
given by Ljunggren and collected in the book \cite[(3.18)]{Gould1972},
a simple proof of which can be found in \cite[Lemma 2.1]{Sun201512}.

\begin{lem}\label{LM:getp}
Let $n\in\mathbb{Z^{+}}$.
Then
\begin{equation}\label{EQ:1}
\sum_{k=0}^{n-1}{k\choose i}{k+j\choose j}
 =\frac{(-1)^jn}{i+j+1}{n-1\choose i}{-n-1\choose j}
\end{equation}
for any $i,j=0,\ldots,n-1$.
\end{lem}
\pf Identity \eqref{EQ:1} can be proved utilizing Gosper's Algorithm in symbolic computation. In fact, by Gosper's algorithm we find that
\[
{k\choose i}{k+j\choose j}
=\Delta_k\left(
\frac{k-i}{i+j+1}{k\choose i}{k+j\choose j}\right),
\]
where $\Delta_k$ denotes the difference operator w.r.t. $k$.
Summing the above identity over $k$ from $0$ to $n-1$, we obtain that
\begin{align*}
\sum_{k=0}^{n-1}{k\choose i}{k+j\choose j}
   = \frac{n-i}{i+j+1}{n\choose i}{n+j\choose j}
   = \frac{(-1)^jn}{i+j+1}{n-1\choose i}{-n-1\choose j}.
\end{align*}
\qed\\
\emph{Proof of the Theorem~\ref{TH:Main1}.}
(1) According to Lemma~\ref{LM:NewFormD}, one can see
\begin{align*}
D_k(r_2,y_2)&=\sum_{j=0}^k{k\choose j}{r_2\choose j}y_2^{j}
             =\sum_{j=0}^{r_2}{r_2\choose j}{k\choose j}y_2^{j}\\
            &=\sum_{j=0}^{r_2}{r_2\choose j}{k+j\choose j}y_2^{j}
                 (1-y_2)^{r_2-j}.
\end{align*}
Then we have
\begin{align}\label{EQ:DD}
& \sum_{k=0}^{p-1}D_k(x_1,y_1)D_k(x_2,y_2)\nonumber\\
\equiv & \sum_{k=0}^{p-1}\sum_{i=0}^k{k\choose i}{r_1\choose i}y_1^{i}
         \sum_{j=0}^{r_2}{r_2\choose j}
                         {k+j\choose j}y_2^{j}(1-y_2)^{r_2-j}\nonumber\\
\equiv & \sum_{k=0}^{p-1}\sum_{i=0}^{p-1}{k\choose i}
               {r_1\choose i}y_1^{i}
         \sum_{j=0}^{p-1}{r_2\choose j}
                         {k+j\choose j}y_2^{j}(1-y_2)^{r_2-j}\nonumber\\
=&\sum_{i=0}^{p-1}\sum_{j=0}^{p-1}{r_1\choose i} {r_2\choose j}
                          y_1^{i}y_2^{j}(1-y_2)^{r_2-j}
   \sum_{k=0}^{p-1}{k\choose i}{k+j\choose j}\nonumber\\
=& \sum_{i=0}^{p-1}\sum_{j=0}^{p-1}{r_1\choose i} {r_2\choose j}
                          y_1^{i}y_2^{j}(1-y_2)^{r_2-j}
         \frac{(-1)^jp}{i+j+1}{p-1\choose i}{-p-1\choose j}\nonumber\\
\equiv & (-1)^{r_2}\sum_{i+j=p-1}{r_1\choose i} {r_2\choose j}
                                 y_1^{i}y_2^{j}(y_2-1)^{r_2-j}\pmod p
\end{align}
with the help of Lemma~\ref{LM:getp}.

When $r=r_1+r_2<p-1$, we know $i+j=p-1$ implies ${r_1\choose i}{r_2\choose j}=0$, and thus
\[
\sum_{k=0}^{p-1}D_k(x_1,y_1)D_k(x_2,y_2)\equiv 0 \pmod p.
\]

When $r=p-1$, it is easy to see that
${r_1\choose i}{r_2\choose j}\neq0$ only when
$i=r_1$ and $j=r_2$.
Then we arrive at
\[
\sum_{k=0}^{p-1}D_k(x_1,y_1)D_k(x_2,y_2)\equiv  (-y_2)^{r_2}y_1^{r_1} \pmod p.
\]
This completes the proof of identity \eqref{EQ:main1}.

(2) To prove identity \eqref{EQ:main2}, we first make the following transformation.
\begin{align*}\label{EQ:DD}
       & \sum_{k=0}^{p-1}S_k(x_1,y_1)S_k(x_2,y_2)\\
\equiv &\sum_{k=0}^{p-1}
\sum_{i=0}^{k}\binom{k}{i}\binom{r_1}{i}\binom{p-1-r_1}{i}y_1^{i}
\sum_{j=0}^{k}\binom{k}{j}\binom{r_2}{j}\binom{p-1-r_2}{j}y_2^{j}\\
\equiv &
\sum_{i=0}^{p-1}\sum_{j=0}^{p-1}
\binom{r_1}{i}\binom{p-1-r_1}{i}\binom{r_2}{j}\binom{p-1-r_2}{j}y_1^{i}y_2^{j}
\sum_{k=0}^{p-1}\binom{k}{i}\binom{k}{j}\pmod p
\end{align*}
If $i>\frac{p-1}{2}$ ($j>\frac{p-1}{2}$), apparently
$\binom{r_1}{i}\binom{p-1-r_1}{i}=0$
($\binom{r_2}{j}\binom{p-1-r_2}{j}=0$).
If $i+j<p-1$, then $\binom{k}{i}\binom{k}{j}$ is a polynomial in $k$ of degree smaller than $p-1$.
This implies that
\[
\sum_{k=0}^{p-1}\binom{k}{i}\binom{k}{j}\equiv 0 \pmod p,
\]
as $\sum_{k=0}^{p-1}k^s\equiv 0 \pmod p$ for any $s=0,1,\ldots,p-2$.
Therefore we only need to consider the case when $i=j=(p-1)/2$, and then
\begin{align}\label{EQ:SS}
\sum_{k=0}^{p-1}S_k(x_1,y_1)S_k(x_2,y_2)
\equiv
c\left(\frac{y_1y_2}{p}\right)
\sum_{k=0}^{p-1}\binom{k}{\frac{p-1}{2}}\binom{k}{\frac{p-1}{2}}\! \pmod p,\!
\end{align}
where $c=\binom{r_1}{\frac{p-1}{2}}\binom{p-1-r_1}{\frac{p-1}{2}}
\binom{r_2}{\frac{p-1}{2}}\binom{p-1-r_2}{\frac{p-1}{2}}$.

When $r_1\neq (p-1)/2$ or $r_2\neq (p-1)/2$, clearly
$c=0$.
This means
\[
\sum_{k=0}^{p-1}S_k(x_1,y_1)S_k(x_2,y_2)\equiv 0 \pmod p
\]
Suppose $r_1=r_2= (p-1)/2$. Then equality \eqref{EQ:SS} leads to
\begin{align*}
\sum_{k=0}^{p-1}S_k(x_1,y_1)S_k(x_2,y_2)
\equiv
\left(\frac{y_1y_2}{p}\right)
\frac{\sum\limits_{k=0}^{p-1}k^{p-1}}{(p-1)!}\binom{p-1}{\frac{p-1}{2}}
\equiv \left(\frac{-y_1y_2}{p}\right)\!\! \pmod p.
\end{align*}
\qed

It is easy to check that
\[
D_n(-\frac{1}{2},-\frac{4}{m})=\sum_{k=0}^{n}\binom{n}{k}\frac{\binom{2k}{k}}{m^k}.
\]
Let $x_1=x_2=-\frac{1}{2}$, $y_1=-\frac{4}{a}$ and $y_2=-\frac{4}{b}$ in theorem \ref{TH:Main1}.
Then $r_1=r_2=\frac{p-1}{2}$ and identity \eqref{EQ:main1} leads to the following corollary.

\begin{coro}
Suppose $p$ is an odd prime, $a,b\in\mathbb{Z}$ and~$\gcd(a,p)=\gcd(b,p)=1$. Then
\begin{equation}\label{EQ:a,b}
\sum_{n=0}^{p-1}\left(\sum_{k=0}^n {n \choose k}\frac{{2k \choose k}}{a^k}\right)
                \left(\sum_{k=0}^n {n \choose k}\frac{{2k \choose k}}{b^k}\right)
\equiv \left(\frac{-ab}{p}\right)
\pmod {p},
\end{equation}
where $\left(\frac{\bullet}{p}\right)$ denotes the Legendre symbol.
\end{coro}
\begin{rem}
The conjectured identity \eqref{CONJ:sun6.7} can be confirmed by taking $a=2$ and $b=-6$ in equality~\eqref{EQ:a,b}.
\end{rem}
\begin{coro}
Let $p$ be a prime, and $p\equiv 1 \pmod 3$. We have
\[
\sum_{n=0}^{p-1}\left(\sum_{k=0}^n {n\choose k}{-1/2 \choose k}a^k \right)
                \left(\sum_{k=0}^n {n\choose k}{-1/3 \choose k}b^k \right)
\equiv 0 \pmod p
\]
for any $p$-adic integers $a,b$.

\pf If $p\equiv 1 \pmod 3$, we know $\langle\frac{-1}{2}\rangle_p=\frac{p-1}{2}$ and
$\langle\frac{-1}{3}\rangle_p=\frac{p-1}{3}$.
Then $r=\frac{5(p-1)}{6}<p-1$, the conclusion follows directly from \eqref{EQ:main1}.
\qed
\end{coro}

Note that
\[
\binom{-1/3}{k}\binom{-2/3}{k}=\frac{\binom{2k}{k}\binom{3k}{k}}{27^k}
\text{ and }
\binom{-1/4}{k}\binom{-3/4}{k}=\frac{\binom{2k}{k}\binom{4k}{2k}}{64^k}.
\]
If $p>3$, then by identity \eqref{EQ:main2} with $x_1=-1/3$ and $x_2=-1/4$ we obtain
\begin{coro}
Let $p$ be an odd prime and $a,b\in\bZ$ are coprime to $p$.
Then
\[
\sum_{k=0}^{p-1}
\left(\sum_{k=0}^{n}
\binom{n}{k}
\frac{\binom{2k}{k}\binom{3k}{k}}{a^k}\right)
\left(\sum_{k=0}^{n}
\binom{n}{k}
\frac{\binom{2k}{k}\binom{4k}{2k}}{b^k}\right)
\equiv 0 \pmod p.
\]
\end{coro}

\section{Proof of Theorem~\ref{TH:Main3}}

To prove Theroem \ref{TH:Main3}, we first recall some known results.
Let $a,m\in\bZ$ with $a>0$ and $\gcd(p,m)=1$, Sun~\cite{Sun2010a} determined $\sum_{k=0}^{p^a-1}{2k\choose k+d}/m^k\pmod {p^2}$ for $d=0,1$ and proved that
\begin{align}
\sum_{k=1}^{p^a-1}\frac{C_k}{m^k}
\equiv
&\ m^{p-1}-1-\frac{m-4}{2}\left(\left(\frac{\Delta}{p^a}\right)-1\right)\nonumber\\
&-\frac{m-4}{2}\left(\frac{\Delta}{p^a-1}\right)u_{p-\left(\frac{\Delta}{p}\right)}(m-2,1) \pmod {p^2},
\end{align}
where $\Delta=m(m-4)$ and $u_n(A,B)$ is the Lucas sequences.

Taking $a=1$, it is straightforward to obtain the following result.
\begin{lem}\label{LM:Catalan}
Let $p$ be an odd prime and $m$ an integer with $\gcd(m,p)=1$. Then we have
\[
\sum_{k=0}^{p-1}\frac{C_k}{m^k}\equiv \frac{4-m}{2}\left(\frac{\Delta}{p}\right)+\frac{m}{2}-1
\pmod {p},
\]
where $C_k$ denotes the $k$th Catalan number and $\Delta=m(m-4)$.
\end{lem}
\noindent
\emph{Proof of Theorem~\ref{TH:Main3}}. Firstly, we do the following decomposition.
\begin{align}
   & \sum_{n=0}^{p-1}\left(\sum_{i=0}^n{n \choose i}\frac{C_i}{a^i}\right)
                    \left(\sum_{j=0}^n{n \choose j}\frac{C_j}{a^j}\right)\nonumber\\
=  & \sum_{i=0}^{p-1}\sum_{j=0}^{p-1}\frac{1}{a^{i+j}(i+1)(j+1)}{2i \choose i}{2j \choose j}\sum_{n=0}^{p-1}{n \choose i}{n \choose j}\nonumber\\
=  & L_1+L_2+L_3, \label{EQ:decomp}
\end{align}
where
\begin{align*}
&L_1=\sum_{i=0}^{p-2}\sum_{j=0}^{p-2}\frac{1}{a^{i+j}(i+1)(j+1)}
                                     {2i \choose i}{2j \choose j}
                               \sum_{n=0}^{p-1}{n \choose i}{n \choose j},\\
&L_2=2\sum_{i=0}^{p-2}\frac{C_i}{p a^{i+p-1}}{2p-2 \choose p-1}
                  {p-1 \choose i},\\
&L_3=\frac{1}{p^2 a^{2p-2}}{2p-2 \choose p-1}^2.
\end{align*}
In the sequel, we will calculate the above three parts separately.
Firstly, notice that
$\binom{2i}{i}\binom{2j}{j}\equiv 0 \pmod p$ if $i>\frac{p-1}{2}$ or $j>\frac{p-1}{2}$,
and that
$
\sum_{k=0}^{p-1}\binom{k}{i}\binom{k}{j}\equiv 0 \pmod p,
$
when $i+j<p-1$.
Then
\[
{2i \choose i}{2j \choose j}\sum_{n=0}^{p-1}{n \choose i}{n \choose j}\not\equiv 0 \pmod p
\]
only when
$i=j=\frac{p-1}{2}$.
Thus we have
\begin{align}\label{EQ:L1}
L_1 & \equiv \frac{4}{(p+1)^2}{p-1 \choose (p-1)/2}^2\cdot
       \sum_{n=0}^{p-1}{n \choose (p-1)/2}^2\nonumber\\
    & \equiv \frac{4\sum_{n=0}^{p-1}n^{p-1}}{(p-1)!}{p-1 \choose (p-1)/2}\nonumber\\
    & \equiv 4\left(\frac{-1}{p}\right) \pmod {p}.
\end{align}
With the help of lemma \ref{LM:Catalan}, we know
\begin{align}\label{EQ:L2}
L_2 & \equiv 2\sum_{i=0}^{p-2}\frac{C_i}{a^i(2p-1)}{2p-1 \choose p}
                     {p-1 \choose i} \nonumber\\
    & \equiv -2\left(\sum_{i=0}^{p-1}\frac{C_i}{(-a)^i}+1\right)\nonumber\\
    & \equiv -(a+4)\left(\frac{a(a+4)}{p}\right)+a \pmod {p}.
\end{align}
Lastly, it is easy to see that
\begin{equation}\label{EQ:L3}
L_3 =\frac{1}{p^2 a^{2p-2}}{2p-2 \choose p-1}^2
      \equiv \frac{1}{(2p-1)^2}{2p-1 \choose p}^2\equiv 1 \pmod {p}.
\end{equation}

Substituting equations~\eqref{EQ:L1}, \eqref{EQ:L2} and \eqref{EQ:L3} back into \eqref{EQ:decomp}, it is straight forward to see that
\[
\sum_{n=0}^{p-1}\left(\sum_{k=0}^n {n \choose k} \frac{C_k}{a^k}\right)^2
\equiv 4\left(\frac{-1}{p}\right)-(a+4)\left(\frac{a(a+4)}{p}\right)+a+1 \pmod {p}.
\]
\qed


\begin{thebibliography}{10}


\bibitem{Gould1972}
H.W. Gould.
\newblock Combinatorial Identities.
\newblock Morgantown Printing and Binding Co., 1972.

\bibitem{GKP1989}
 R.L. Graham, D.E. Knuth and  O. Patashnik.
\newblock Concrete Mathematics.
\newblock Addison-Wesley, New York, 1989.


\bibitem{Guo2016}
V.J.W. Guo.
\newblock Proof of Sun's conjectures on integer-valued polynomals.
\newblock J. Math. Anal. Appl. 444(2016), 182--191.


\bibitem{HouWang2018}
Q.-H. Hou and Y.-S. Wang.
\newblock Constant term evaluation and two kinds of congruence.
\newblock Int. J. Number Theory 14(2018), 2013--2022.

\bibitem{KW2006}
K. Kimoto and M. Wakayama.
\newblock Ap\'ery-like numbers arising from special values of spectral zeta function for non-commutative harmonic oscillators.
\newblock Kyushu J. Math. 60(2006), 383--404.

\bibitem{Liu2016}
J.-C. Liu.
\newblock A generalized supercongruence of Kimoto and Wakayama.
\newblock J. Math. Anal. Appl. 467(2016), 15--25.

\bibitem{Liu2017}
J.-C. Liu.
\newblock Proof of some divisibility results on sums involving binomial coefficients.
\newblock J. Number Theory 180(2017), 566--572.




\bibitem{LOS2014}
L. Long, R. Osburn and H. Swisher.
\newblock On a conjecture of Kimoto and Wakayama.
\newblock Proc. Amer. Math. Soc. 144(2016), 4319--4327.


\bibitem{Mortenson2003a}
E. Mortenson.
\newblock A supercongruence conjecture of Rodriguez-Villegas for a certain truncated hypergeometric function.
\newblock J. Number Theory 99(2003), 139--147.


\bibitem{Mortenson2003b}
 E. Mortenson.
\newblock Supercongruences between truncated $_2F_1$\  by geometric functions and their Gaussian analogs.
\newblock Trans. Amer. Math. Soc. 355(2003), 987--1007.



\bibitem{RV2003}
F. Rodriguez-Villegas.
\newblock Hypergeometric families of Calabi--Yau manifolds.
\newblock Calabi--Yau Varieties and Mirror Symmetry (Toronto, ON, 2001), pp. 223--231, Fields Inst. Commun. Vol. 38, Amer. Math. Soc., Providence, RI, 2003.

%
%
%



\bibitem{Sun2003}
Z.-W. Sun.
\newblock Combinatorial identities in dual sequences.
\newblock Eueopean J. Combin. 24(2003), 709--718.
\bibitem{Sun2010a}
Z.-W. Sun.
\newblock Binomial coefficients, Catalan numbers and Lucas quotients.
\newblock Sci. China Math. 53(2010), 2473--2488.


\bibitem{Sun201512}
Z.-W. Sun.
\newblock Supercongruences involving dual sequences.
\newblock Finite Fields Appl. 46(2017), 179--216.



\end{thebibliography}
\end{document}